\pdfoutput=1
\documentclass[11pt]{amsart}
\usepackage{graphicx, verbatim, fancyhdr,indent} 
\usepackage[nohug]{diagrams}\diagramstyle[labelstyle=\scriptstyle]
\usepackage[usenames,dvipsnames]{color}
\usepackage[colorlinks=true, urlcolor=MidnightBlue, linkcolor=MidnightBlue, citecolor=MidnightBlue, pdftitle={Triple linking numbers, ambiguous Hopf invariants and integral formulas for three-component links}, pdfauthor={Dennis DeTurck, Herman Gluck, Rafal Komendarczyk, Paul Melvin, Clayton Shonkwiler and David Shea Vela-Vick}, pdfsubject={Geometric Topology, Knot Theory}, pdfkeywords={mu invariant, Hopf invariant, link-homotopy, string links, helicity}]{hyperref}

\bfseries
\def\part#1{\subsection{#1}}

\newtheorem{theorem}{Theorem}

\newtheorem*{proposition}{Proposition}

\theoremstyle{definition}
\newtheorem*{remark}{Remark}

\newcommand{\thmref}[1]{Theorem~\ref{#1}}
\newcommand{\secref}[1]{Section~\ref{#1}}

\newcommand{\foot}{\setcounter{footnote}{1}\footnote}

\newcommand{\itb}{\bfseries\itshape}

\def\bz{\mathbb Z}

\def\br{\mathbb R}

\def\boldF{\mathbb F}

\def\ba{\mathbb A}
\def\bl{\mathbb L}

\def\calg{\mathcal G}
\def\calh{\mathcal H}
\def\M{\text{Maps}(S^1\times S^1,S^2)}
\def\Mp{\text{Maps}_p(S^1\times S^1,S^2)}

\def\tor{T^3}
\def\conf{\mathrm{Conf}_3S^3}
\def\gras{\mathrm{G}_2\br^4}
\def\dt{\,\raisebox{.2ex}{\text{\tiny\textbullet}}\,}
\def\dtw{\,\dt\,}

\def\col{\colon}
\def\st{\ | \ }
\def\bx{x^{-1}} % {\bar x}
\def\pr{\mathrm{\wp}} 
\def\g{\widetilde g}

\begin{document}
%%%%%%%%%%%%%%%

%%%%%%%%%%%%%%%%%%%%%%%%%
%% TITLE and ABSTRACT
%%%%%%%%%%%%%%%%%%%%%%%%%

\title{Triple linking numbers, ambiguous Hopf invariants \\ {\footnotesize and} integral formulas {\footnotesize for} three-component links}

\author{Dennis DeTurck \and Herman Gluck \and Rafal Komendarczyk \and \\ Paul Melvin \and Clayton Shonkwiler \and David Shea Vela-Vick}

\dedicatory{To Manfredo do Carmo \\ in friendship and admiration, \!on his $80^{\text{th}}$ birthday}

\date{12 January, 2009}

\begin{abstract}
	
	Three-component links in the $3$-dimensional sphere were classified up to link homotopy by John Milnor in his senior thesis, published in 1954.  A complete set of invariants is given by the pairwise linking numbers $p$, $q$ and $r$ of the components, and by the residue class of one further integer $\mu$, the ``triple linking number'' of the title, which is well-defined modulo the greatest common divisor of $p$, $q$ and $r$.

	\vskip .1in 

	%%% FIGURE: Borromean Rings %%%%%%%%%%%%%%%%%%%%%%%%%% 
	\begin{center} \includegraphics[height=135px]{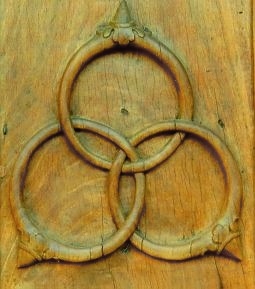} \end{center} %replaced cremona2.jpg
	\begin{center} {\footnotesize The Borromean rings: $p=q=r=0$, $\mu=\pm1$} \end{center}
	%%%%%%%%%%%%%%%%%%%%%%%%%%%%%%%%%%%%%%%%%%%
	\vskip .1in
	
	To each such link $L$ we associate a geometrically natural characteristic map $g_L$ from the $3$-torus to the $2$-sphere in such a way that link homotopies of $L$ become homotopies of $g_L$.  Maps of the 3-torus to the 2-sphere were classified up to homotopy by Lev Pontryagin in 1941.  A complete set of invariants is given by the degrees $p$, $q$ and $r$ of their restrictions to the $2$-dimensional coordinate subtori, and by the residue class of one further integer $\nu$, an ``ambiguous Hopf invariant" which is well-defined modulo twice the greatest common divisor of $p$, $q$ and $r$.
	\vspace{4pt}
	
	We show that the pairwise linking numbers $p$, $q$ and $r$ of the components of $L$ are equal to the degrees of its characteristic map $g_L$ restricted to the $2$-dimensional subtori, and that twice Milnor's $\mu$-invariant for $L$ is equal to Pontryagin's $\nu$-invariant for $g_L$.
	\vspace{4pt}
	
	When $p$, $q$ and $r$ are all zero, the $\mu$- and $\nu$-invariants are ordinary integers.  In this case we use J.\,H.\,C.~Whitehead's integral formula for the Hopf invariant, adapted to maps of the $3$-torus to the $2$-sphere, together with a formula for the fundamental solution of the scalar Laplacian on the $3$-torus as a Fourier series in three variables, to provide an explicit integral formula for $\nu$, and hence for $\mu$.
	\vspace{4pt}
	
	We give here only sketches of the proofs of the main results, with full details to appear elsewhere.
\end{abstract}

\maketitle

\parskip 4pt

\topmargin=16pt
\textheight=560pt
\headheight=14pt 

%%%%%%%%%%%%%%%%%%%%%%%%%
%% SECTION 1: Statement of Results
%%%%%%%%%%%%%%%%%%%%%%%%%

\section{Statement of results}
\label{sec:results}

Consider the configuration space 
\[
\conf \subset S^3 \times S^3 \times S^3
\]
of ordered triples $(x,y,z)$ of distinct points in the unit $3$-sphere $S^3$ in $\br^4$.  Since $x$, $y$ and $z$ are distinct, they span a $2$-plane in $\br^4$.  Orient this plane so that the vectors from $x$ to $y$ and from $x$ to $z$ form a positive basis, and then move it parallel to itself until it passes through the origin.  The result is an element $G(x,y,z)$ of the Grassmann manifold $\gras$ of all oriented 2-planes through the origin in $\br^4$.  This defines the {\itb Grassmann map} 
$$
G: \conf \longrightarrow \gras.
$$
It is equivariant with respect to the diagonal $O(4)$ action on $S^3 \times S^3 \times S^3$ and the usual $O(4)$ action on $\gras$.

\vspace{.2in}

%%% FIGURE 1: Grassmann Map %%%
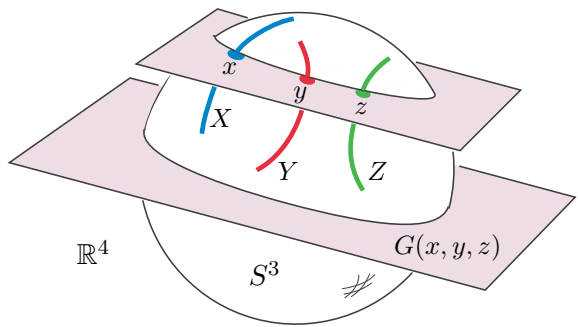
\begin{figure}[h!]
\includegraphics[height=120pt]{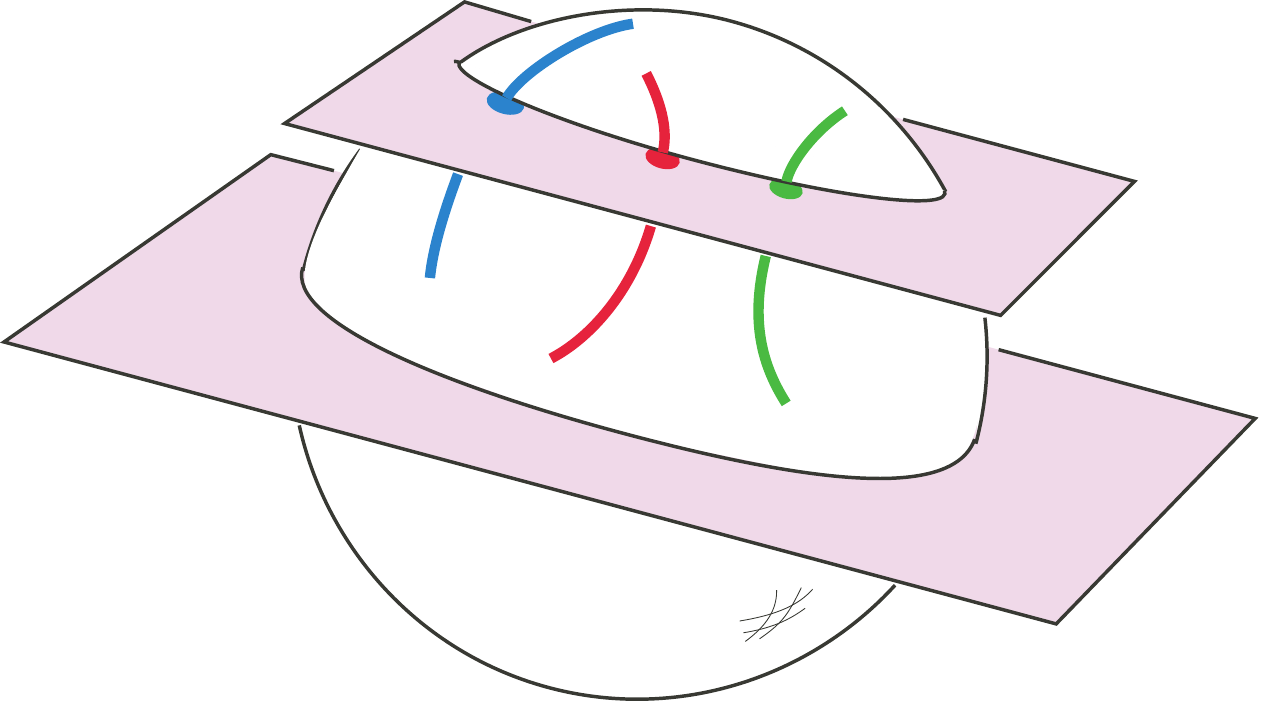}
\put(-70,27){\small $G(x,y,z)$}
\put(-135,95){\small $x$}
\put(-108,86){\small $y$}
\put(-85,80){\small $z$}
\put(-140,75){\small $X$}
\put(-114,55){\small $Y$}
\put(-80,55){\small $Z$}
\put(-125,15){$S^3$}
\put(-190,23){$\br^4$}
\caption{The Grassmann map}
\label{fig:gras}
\end{figure}
%%%%%%%%%%%%%%%%%%%%

The Grassmann manifold $\gras$ is isometric (up to scale) to the product $S^2 \times S^2$ of two unit 2-spheres.  Let $\pi\colon\gras\to S^2$ denote orthogonal projection to either factor.   

Given any ordered oriented link $L$ in $S^3$ with three parametrized components
\[
	X = \{x(s) \st s \in S^1\} \ , \ Y = \{y(t) \st t \in S^1 \} \ \ \text{and}\ \  Z = \{z(u) \st u \in S^1 \},
\]
where $S^1$ is the unit circle in $\mathbb{R}^2$, we define the {\itb characteristic map} of $L$
\[
	g_L\colon T^3 = S^1 \times S^1 \times S^1 \longrightarrow S^2
\]
by $g_L(s,t,u) = \pi(G(x(s),y(t), z(u)))$.  In \secref{sec:characteristicmap} we give an explicit formula for this map as the unit normalization of a vector field on $\tor$ whose components are quadratic polynomials in the components of $x(s)$, $y(t)$ and $z(u)$.  

The homotopy class of $g_L$ is unchanged under  any {\itb link homotopy} of $L$, meaning a deformation during which each component may cross itself, but different components may not intersect.  

\pagebreak

\begin{theorem}\label{thm:A}
	\textbf{The pairwise linking numbers \boldmath{$p$}, $q$ and $r$ of the link $L$ are equal to the degrees of its characteristic map $g_L$ on the 2-dimensional \mbox{coordinate} subtori of $T^3$, while twice Milnor's $\mu$-invariant for $L$ is equal to \mbox{Pontryagin's} $\nu$-invariant for $g_L$.}
\end{theorem}
\vspace{.01in}

\begin{remark}
	Milnor's $\mu$-invariant, typically denoted $\overline{\mu}_{123}$, is descriptive of a {\it single} three-component link.  In contrast, Pontryagin's $\nu$-invariant is the cohomology class of a difference cocycle comparing {\it two} maps from $T^3$ to $S^2$ that are homotopic on the $2\text{-skeleton}$ of $T^3$.  In particular, it assigns to any pair $g$, $g'$ of such maps whose degrees on the coordinate 2-tori are $p$, $q$ and $r$, an integer $\nu(g,g')$ that is well-defined modulo $2\gcd(p,q,r)$.  With this understanding, the last statement in \thmref{thm:A} asserts that
	\[
		2(\mu(L) - \mu(L')) \ \equiv \ \nu(g_L,g_{L'}) \mod{2\gcd(p,q,r)},
	\]
for any two links $L$ and $L'$ with the same pairwise linking numbers $p$, $q$ and $r$.
\end{remark}
\vspace{.05in}

We will sketch here two quite different proofs of \thmref{thm:A}, a topological one in \secref{sec:top_proof} using framed cobordism of framed links in the 3-torus, and an algebraic one in \secref{sec:alg_proof} using the group of link homotopy classes of three-component string links and the fundamental groups of spaces of maps of the 2-torus to the 2-sphere.

To state the integral formula for Milnor's $\mu$-invariant when the pairwise linking numbers are zero, let $\omega$ denote the Euclidean area 2-form on $S^2$, normalized so that the total area is 1 instead of $4\pi$.  Then $\omega$ pulls back under the characteristic map $g_L$ to a closed 2-form on $T^3$, which can be converted in the usual way to a divergence-free vector field $V_L$ on $T^3$.  In \secref{sec:thmB} we give  explicit formulas for $V_L$, and also for the fundamental solution $\varphi$ of the scalar Laplacian on the 3-torus as a Fourier series in three variables.  These are the key ingredients in the integral formula below.

\vspace{.05in}
\begin{theorem}\label{thm:B}
	\textbf{\boldmath{If the pairwise linking numbers $p$, $q$ and $r$ of the three components of $L$ are all zero, then Milnor's $\mu$-invariant of $L$ is given by the formula
			\[
				\mu(L) \ = \  \frac12\, \int_{T^3 \times T^3} V_L(\sigma) \times V_L(\tau) \dtw \nabla_\sigma		         \varphi\left(\sigma - \tau\right)\ d\sigma\, d\tau .
			\]}}
\end{theorem}
\vspace{.05in}

Here $\nabla_\sigma$ indicates the gradient with respect to $\sigma$, the difference $\sigma - \tau$ is taken in the abelian group structure of the torus, and $d\sigma$ and $d\tau$ are volume elements.  The integrand is invariant under the action of the group $SO(4)$ of orientation-preserving rigid motions of $S^3$ on the link $L$, attesting to the naturality of the formula.  We will see in the next section that the integral above expresses the ``helicity'' of the vector field $V_L$ on $T^3$.

\clearpage

%%%%%%%%%%%%%%%%%%%%%%%%%
%% SECTION 2: Background and Motivation
%%%%%%%%%%%%%%%%%%%%%%%%%

\section{Background and motivation}
\label{sec:background}

Let $L$ be an ordered oriented link in $\br^3$ with two parametrized components 
\[
	X = \{x(s) \st s \in S^1\} \ \  \text{and} \ \ Y = \{y(t) \st t \in S^1 \}.
\]
The classical linking number $\text{Lk}(X, Y)$ is the degree of the Gauss map $S^1\times S^1\to S^2$ sending $(s,t)$ to $(y(t)-x(s))/\|y(t)-x(s)\|$, and can be expressed by the famous integral formula of Gauss~\cite{Gauss},
\begin{align*}
	\text{Lk}(X, & Y) = 
	\frac{1}{4\pi} \int_{S^1 \times S^1} x'(s)  \times  y'(t)  \dtw   \frac{x(s)-y(t)}{\|x(s)-y(t)\|^3} \ ds \, dt \\
	& = \int_{S^1 \times S^1} x'(s)  \times  y'(t)  \dtw  \nabla_x \, \varphi\left(\|x(s)-y(t)\|\right) \ ds \, dt,
\end{align*}
where $\varphi(r) = -1/(4\pi r)$ is the fundamental solution of the scalar Laplacian in $\br^3$.  The integrand is invariant under the group of orientation-preserving rigid motions of $\br^3$, acting on the link $L$.  Corresponding formulas in $S^3$ appear in DeTurck and Gluck~\cite{DeTurckGluck} and in Kuperberg~\cite{Kuperberg}.

Theorems~\ref{thm:A} and~\ref{thm:B} above give a similar formulation of Milnor's triple linking number in $S^3$.  We emphasize that these two theorems are set specifically in $S^3$, and that so far we have been unable to find corresponding formulas in Euclidean space $\br^3$ which are equivariant (for \thmref{thm:A}) and invariant (for \thmref{thm:B}) under the {\sl noncompact} group of orientation-preserving rigid motions of $\br^3$.  

For some background on higher order linking invariants, see Milnor \cite{Milnor57} and, for example, Massey~\cite{Massey68}, Casson~\cite{Casson75}, Turaev~\cite{Turaev76}, Porter~\cite{Porter}, Fenn~\cite{Fenn83}, Orr~\cite{Orr}, and Cochran~\cite{Cochran}.

% Also Levine?  Not sure about Rozansky (1994), Thurston (1995) and Akhmetiev (1998 and 2005).  

The {\itb helicity} of a vector field $V$ defined on a bounded domain $\Omega$ in $\br^3$ is given by the formula
\[
	\begin{aligned}
	\text{Hel}&(V) = \int_{\Omega \times \Omega} V(x) \times V(y) \dtw \frac{x-y}{\|x-y\|^3} \  dx \, dy \\
	&= \int_{\Omega \times \Omega} V(x) \times V(y) \dtw \nabla_x \varphi\left(\|x-y\|\right) \ dx \, dy
	\end{aligned}
\]

\noindent where, as above, $\varphi$ is the fundamental solution of the scalar Laplacian on $\br^3$.  

Woltjer~\cite{Woltjer} introduced this notion during his study of the magnetic field in the Crab Nebula, and showed that the helicity of a magnetic field remains constant as the field evolves according to the equations of ideal magnetohydrodynamics, and that it provides a lower bound for the field energy during such evolution.  The term ``helicity'' was coined by Moffatt~\cite{Moffatt}, who also derived the above formula.  

There is no mistaking the analogy with Gauss's linking integral, and no surprise that helicity is a measure of the extent to which the orbits of $V$ wrap and coil around one another.  Since its introduction, helicity has played an important role in astrophysics and solar physics, and in plasma physics here on earth.

Looking back at \thmref{thm:B}, we see that the integral in our formula for Milnor's $\mu$-invariant of a three-component link $L$ in the 3-sphere expresses the helicity of the associated vector field $V_L$ on the 3-torus.

Our study was motivated by a problem proposed by Arnol$'$d and Khesin~\cite{ArnoldKhesin98} regarding the search for ``higher helicities'' for divergence-free vector fields.  In their own words:
\vspace{.05in}
\begin{indentation}{2.6em}{2.6em}
\small\itb
\noindent The dream is to define such a hierarchy of invariants for generic vector fields such that, whereas all the invariants of order $\leq k$ have zero value for a given field and there exists a nonzero invariant of order $k+1$, this nonzero invariant provides a lower bound for the field energy.
\end{indentation}
\vspace{.05in}

Many others have been motivated by this problem, and have contributed to its understanding; see, for example, Berger and Field~\cite{Berger84}, Berger~\cite{Berger90, Berger91}, Evans and Berger~\cite{EvansBerger92}, Akhmetiev and Ruzmaiken~\cite{Akhmetiev94, Akhmetiev95}, Akhmetiev~\cite{Akhmetiev98}, Laurence and Stredulinsky~\cite{Laurence}, Hornig and Mayer~\cite{Hornig}, Rivi\`ere~\cite{Riviere}, Khesin~\cite{Khesin}, and Bodecker and Hornig~\cite{Bodecker04}.

% Also Rozansky (1994), Auckly et al?  

The formulation in Theorems~\ref{thm:A} and~\ref{thm:B} has led to partial results that address the case of vector fields on invariant domains such as flux tubes modeled on the Borromean rings; see Komendarczyk~\cite{Komendarczyk}.

\vspace{.2in}

\subsection*{Acknowledgements} We are grateful to Fred Cohen and Jim Stasheff for their substantial input and help during the preparation of this paper, and to Toshitake Kohno, whose 2002 paper provided the original inspiration for this work.  Komendarczyk and Shonkwiler also acknowledge support from DARPA grant \#FA9550-08-1-0386.

The Borromean rings shown on the first page of this paper are from a panel in the carved walnut doors of the Church of San Sigismondo in Cremona, Italy.  The photograph is courtesy of Peter Cromwell.

\clearpage

%%%%%%%%%%%%%%%%%%%%%%%%%
%% SECTION 3: Explicit formula for the characteristic map
%%%%%%%%%%%%%%%%%%%%%%%%%

\section{Explicit formula for the characteristic map $g_L$} 
\label{sec:characteristicmap}

View $\br^4$ as the space of quaternions, with $1$, $i$, $j$, $k$ as basis, and consider the vector field $F: \conf \to \br^3 - \{0\}$ defined by
\[
F(x,y,z) \ = \ \begin{pmatrix} x \dt iy + y \dt iz + z \dt ix \\ x \dt jy + y \dt jz + z \dt jx \\ x\dt ky + y \dt kz + z \dt kx \end{pmatrix}.
\]
Here $\dt$ denotes the dot product in $\br^4$.  The components of $F(x,y,z)$ are quadratic polynomials in the components of $x$, $y$ and $z$ in $\br^4$, and the norm $\|F(x,y,z)\|$ is twice the area of the triangle in $\br^4$ with vertices at $x$, $y$ and $z$. 

Now let $L$ be a three-component link in $S^3$ with parametrized components 
\[
	X = \{x(s) \st s \in S^1\} \ , \ Y = \{y(t) \st t \in S^1 \} \ \ \text{and}\ \ Z = \{z(u) \st u \in S^1 \}.
\]
This defines an embedding \ $e_L:\tor\hookrightarrow\conf$ \ given by $e_L(s,t,u) = (x(s),y(t),z(u))$.  We compute that the characteristic map $g_L\col\tor \to S^2$ is the unit normalization of the composition $F\circ e_L$, that is,
\[
	g_L(s,t,u) = \frac{F(x(s),y(t),z(u))}{\|F(x(s),y(t),z(u))\|}.
\]

Note that the $g_L$ is ``symmetric" in the sense that  it transforms under any permutation of the components of $L$ by precomposing with the corresponding permutation automorphism of $\tor$, and then multiplying by the sign of the permutation.

%%%%%%%%%%%%%%%%%%%%%%%%%%%%
\part{An asymmetric version of the characteristic map} 
\label{sec:alt_char_map}
%%%%%%%%%%%%%%%%%%%%%%%%%%%%

\begin{figure}[h!]
\includegraphics[height=150pt]{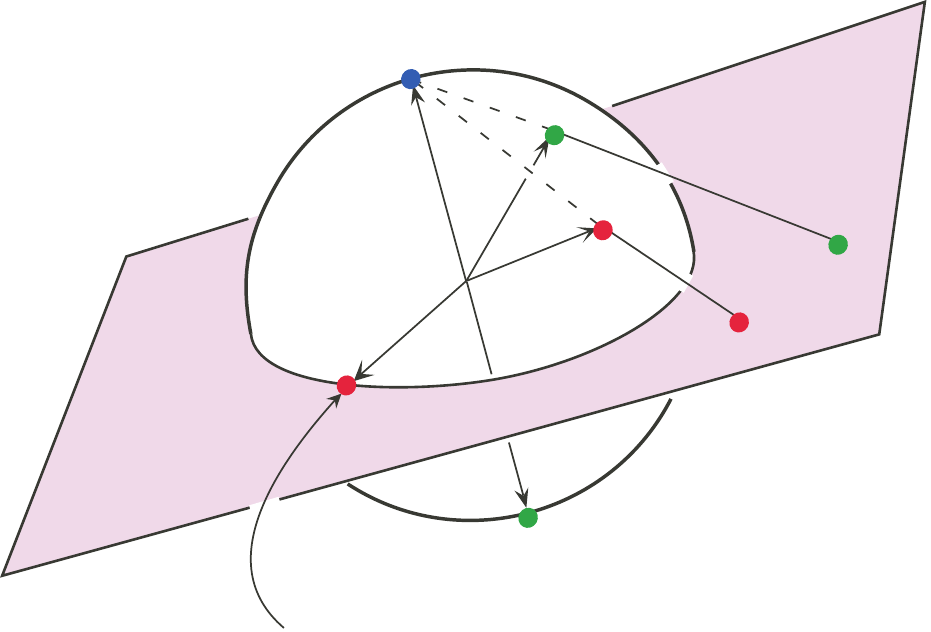}
\put(-127,136){\small$x$}
\put(-105,18){\small$-x$}
\put(-80,103){\small$y$}
\put(-92,123){\small$z$}
\put(-40,72){\small $\pr_xy$}
\put(-30,100){\small $\pr_xz$}
\put(-187,45){$\br^3_x$}
\put(-145,90){$S^3$}
\put(-150,0){\small$(\pr_xy-\pr_xz)/\| \pr_xy-\pr_xz\|$}
\caption{Stereographic projection}
\label{fig:retr}
\end{figure}
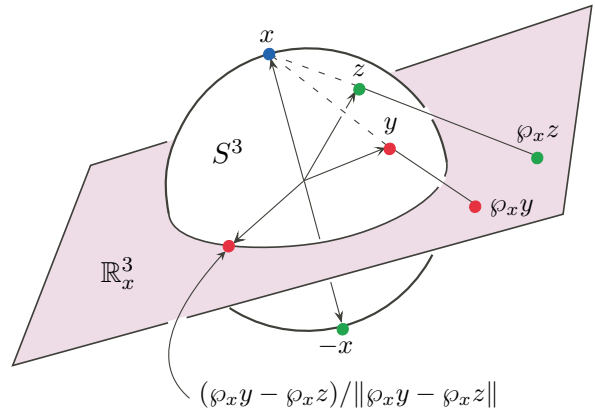

Continuing to view \mbox{$S^3 \subset \br^4$}, let $\pr_x$ denote stereographic projection of $S^3 - \{x\}$ onto the 3-space $\br_x^3$ through the origin in $\br^4$ that is orthogonal to the vector $x$, as shown in Figure~\ref{fig:retr}.  
%%% FIGURE 2: Retraction %%%
%%%%%%%%%%%%%%%%%

For any $(x,y,z)\in\conf$, consider the points $\pr_xy$ and $\pr_xz$ in $\br^3_x$.  Translation in $\br^3_x$ moves $\pr_xz$ to the origin, and then dilation in $\br^3_x$ makes the translated $\pr_xy$ into a unit vector.  Composing with $\pr_x^{-1}$, we see that $x$ has been kept fixed, $y$ has moved to the point $(\pr_xy-\pr_xz)/\| \pr_xy-\pr_xz\|$ on the equatorial $2$-sphere $S^2_x = S^3 \cap \br^3_x$, and $z$ has moved to $-x$, as indicated in the figure.  

This procedure defines a deformation retraction
$$
r(x,y,z) \ = \ (x \, , \, (\pr_xy-\pr_xz)/\| \pr_xy-\pr_xz\|\ \, ,\, -x)
$$
of $\conf$ onto the subspace $\{(x,w,-x) \st x \dt w = 0\}$, which is a copy of the unit tangent bundle $US^3$ of the $3$-sphere via the correspondence $(x,w,-x) \leftrightarrow (x,w)$.  Let $\pi:US^3\to S^2$ denote the projection onto the fiber, sending $(x,w)$ to $w\bx$.
%where $\bx$ is the quaternionic conjugate of $x\in S^3$

Now define the {\itb asymmetric characteristic map} \ $\g_L: \tor\longrightarrow S^2$ \ of a link $L$, as above, to be the composition $\pi\circ r\circ e_L$.  Noting that $(\pr_xv)\bx = \pr_1(v\bx)$, we have explicitly
\[
\g_L(s,t,u) \ = \ \frac{\pr_1(y\bx)-\pr_1(z\bx)}{\|\pr_1(y\bx)-\pr_1(z\bx)\|}
\]

\noindent where $x = x(s)$, $y = y(t)$ and $z = z(u)$ parametrize the components of $L$.   This map is easily seen to be homotopic to the characteristic map $g_L$ defined above.   

The restriction of $\g_L$ to $*\times S^1\times S^1$ is the negative of the Gauss map for the link $\pr_1((Y\cup Z)\bx)$, and so noting that $\pr_1$ is orientation reversing (with the usual sign conventions) its degree is the linking number of $Y$ with $Z$.  Since  $g_L$  is homotopic to $\g_L$,  the same is true for $g_L$.  But then  it follows from the symmetry of $g_L$ that its degree on $S^1\times * \times S^1$ is the linking number of $X$ with $Z$,  and its degree on $S^1 \times S^1 \times *$  is the linking number of $X$ with $Y$.  This proves the first statement in  \thmref{thm:A}.

This version $\g_L$ of the characteristic map will also facilitate the topological proof of the rest of \thmref{thm:A}, to be given next.

\clearpage

%%%%%%%%%%%%%%%%%%%%%%%%%
%% SECTION 4: Sketch of the topological proof of Theorem A
%%%%%%%%%%%%%%%%%%%%%%%%%

%%%%%%%%%%%%%%%%%%%%%%%%%%%%%%
\section{Sketch of the topological proof of \thmref{thm:A}}\label{sec:top_proof}
%%%%%%%%%%%%%%%%%%%%%%%%%%%%%%

\begin{figure}[h!]
\includegraphics[height=100pt]{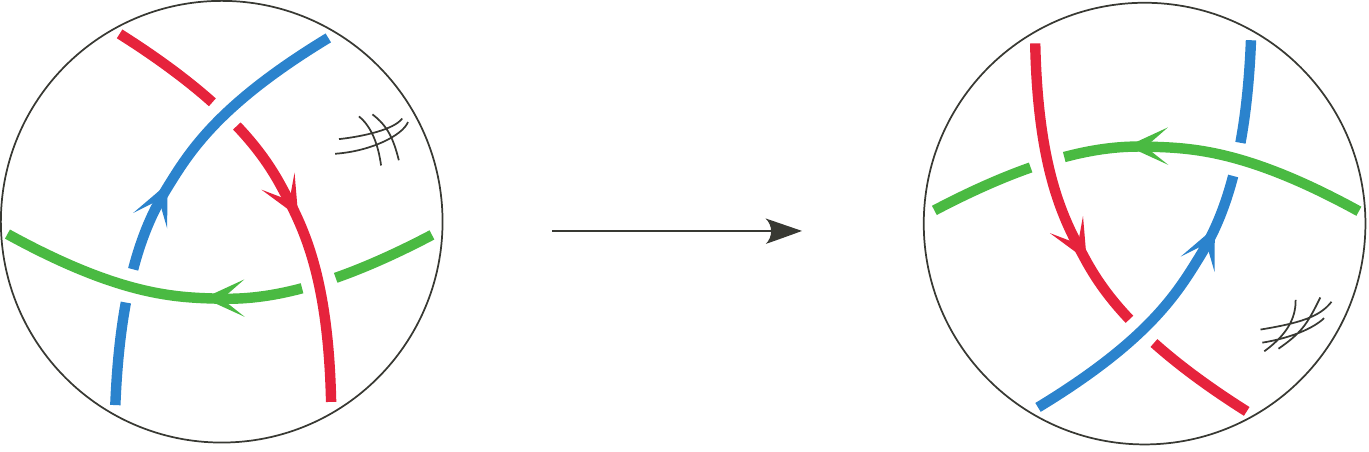}
\put(-228,95){$X$}
\put(-28,99){$X'$}
\put(-230,-5){$Y$}
\put(-20,-3){$Y'$}
\put(-320,45){$Z$}
\put(-115,46){$Z'$}
%\put(-250,9){$\small B$}
%\put(-87,9){$\small B$}
\caption{The delta move $L\to L'$}
\label{fig:delta}
\end{figure}

Starting with a link in the $3$-sphere, consider the {\itb delta move} shown in Figure \ref{fig:delta}, which may be thought of as a higher order variant of a crossing change.  It takes place within a $3$-ball, outside of which the link is left fixed.

%%% FIGURE 3: Delta Move %%%
%%%%%%%%%%%%%%%%%% 

The delta move was introduced by Murakami and Nakanishi~\cite{MurakamiNakanishi}, who proved that an appropriate sequence of such moves can transform any link into any other link with the same number of components, provided the two links have the same pairwise linking numbers.

The key organizational idea for this proof of \thmref{thm:A} is to show that the delta move, when applied to the three components $X$, $Y$ and $Z$ of the link $L$, as shown in Figure~\ref{fig:delta}, increases its Milnor $\mu$-invariant by $1$, while increasing the Pontryagin $\nu$-invariant of its characteristic map $g_L$ by $2$.

\begin{figure}[h!]
\includegraphics[height=120pt]{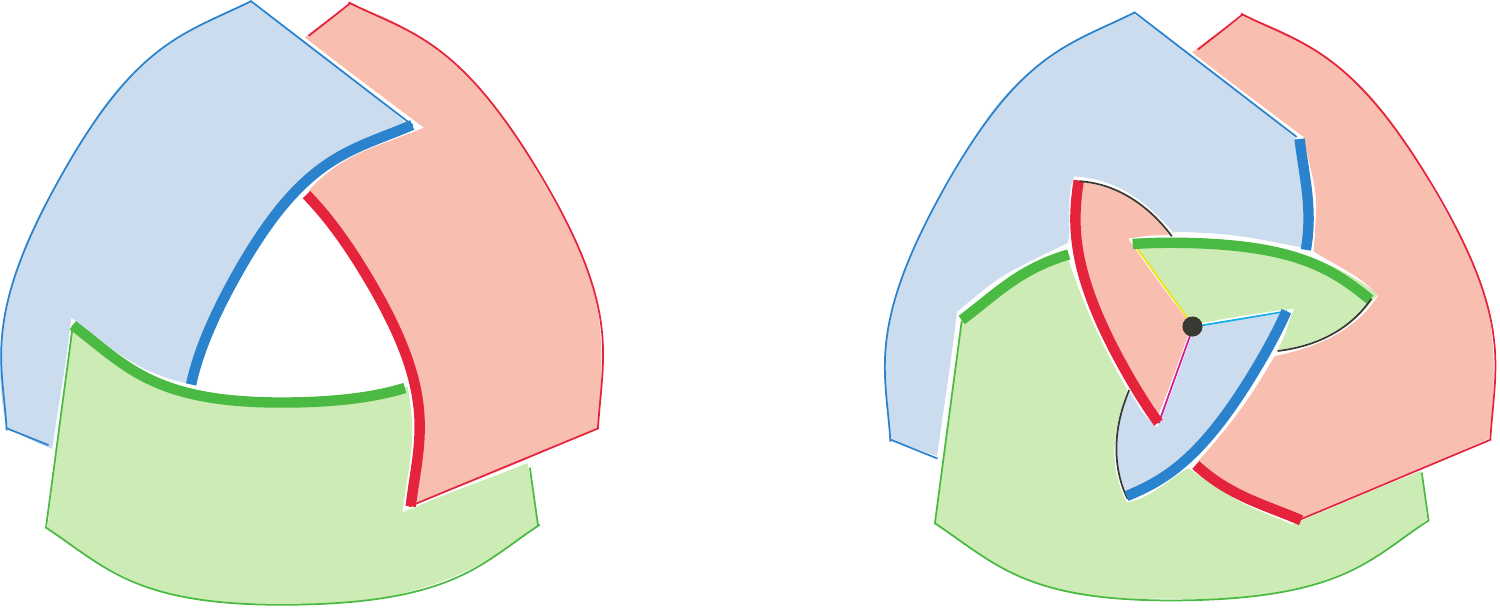}
\caption{Seifert surfaces for $L$ and $L'$}
\label{fig:seifert}
\end{figure}

The fact that the delta move increases $\mu$ by $1$ is well known to experts; our proof relies on the geometric formula for $\mu$ due to Mellor and Melvin~\cite{MellorMelvin} in terms of how each link component pierces the Seifert surfaces of the other two components, plus a count of the triple point intersections of these surfaces.  In particular, one can use a family of Seifert surfaces that differ only near the delta move, as shown in Figure~\ref{fig:seifert}.

%%% FIGURE 4: Seifert Surfaces %%%
%%%%%%%%%%%%%%%%%%%%
 
To see how the delta move affects the $\nu$-invariant, we will view $\nu$ as a relative Euler class, following Gompf~\cite{Gompf} and Cencelj, Repov{\v s} and Skopenkov~\cite{Cencelj}.  

To that end let $L$ and $L'$ be two three-component links in $S^3$ with the same pairwise linking numbers $p$, $q$ and $r$, and let $\bl$ and $\bl'$ be the framed links in the $3$-torus\foot{For the reader's convenience, all subsets of $\tor$ in this section are written in blackboard bold to distinguish them from subsets of $S^3$.} that are the preimages of a common regular value of their characteristic maps $g_L$ and $g_{L'}$.  Orient $\bl$ and $\bl'$ so that, when combined with the pullback of the orientation on $S^2$ to a tangent $2$-plane transverse to these links, we get the given orientation on $\tor$.  Since $L$ and $L'$ have the same pairwise linking numbers, $\bl$ and $\bl'$ are homologous in $\tor$.  It then follows by a standard argument that $\bl\times0$ and $\bl'\times1$ together bound an embedded surface $\boldF$ in $\tor\times[0,1]$.   

The {\itb relative normal Euler class} $e(\boldF)$ is the intersection number of $\boldF$ with a generic perturbation of itself that is directed by the given framings along $\bl\times0$ and $\bl'\times1$, but is otherwise arbitrary.  Then, according to Gompf and Cencelj--Repov{\v s}--Skopenkov, 
$$
\nu(g_{L'},g_L) \ \equiv \ e(\boldF) \ \mod{2\gcd(p,q,r)}.
$$

The key step in seeing how the delta move affects the $\nu$-invariant is to adjust $L$ and $L'$ by link homotopies so that up to isotopy
$$
\bl' = \bl \cup \bl^*
$$
where $\bl^*$ is a two-component link bounding an annulus $\ba \subset \tor - \bl$ with relative normal Euler class 2.   Then the surface $\boldF \ = \ (\bl\times[0,1])\cup(\ba\times 1)$, with $\partial\boldF = \bl'\times1-\bl\times0$, has $e(\boldF)=2$, since $e(\bl\times[0,1])=0$.  Thus $\nu(g_{L'},g_L) \equiv e(\boldF) = 2$, showing that the delta move increases the $\nu$-invariant of $L$ by $2$.  Explaining how this step is carried out will complete our sketch of the topological proof of \thmref{thm:A}. 

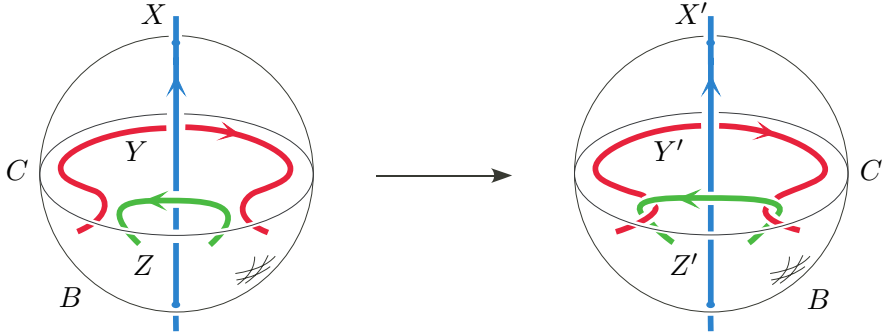
\begin{figure}[h!]
\includegraphics[height=120pt]{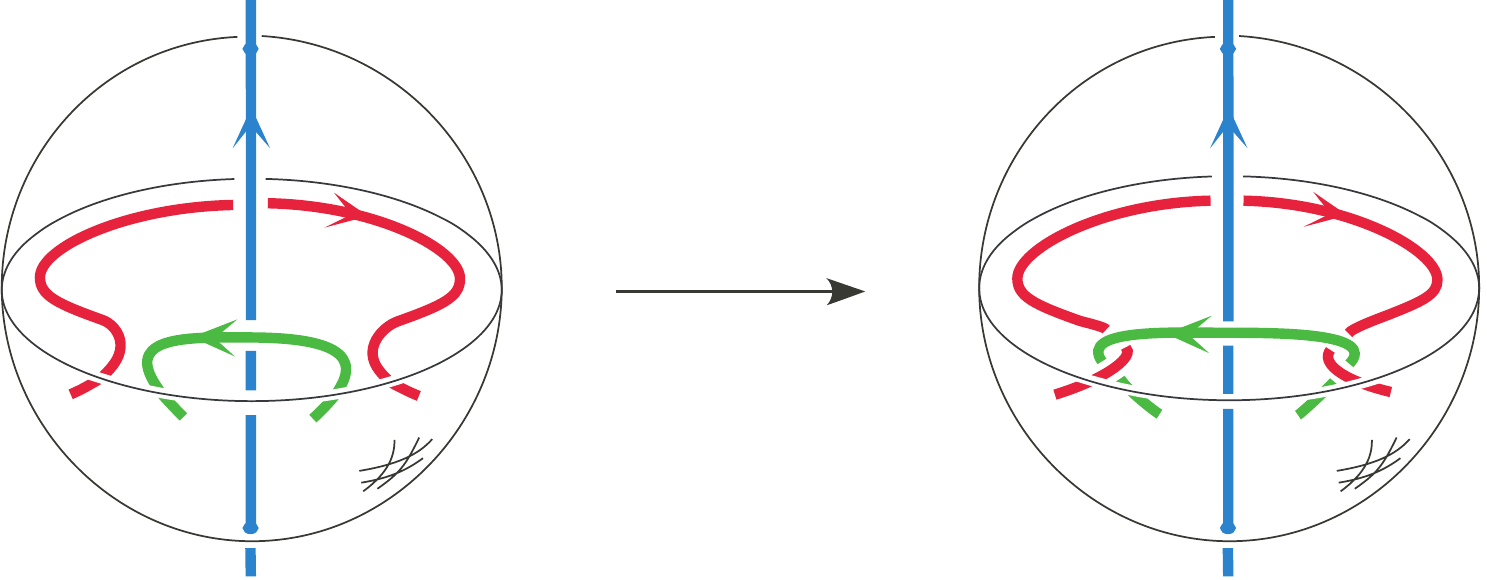}
\put(-300,10){$B$}
\put(-17,9){$B$}
\put(-320,58){$C$}
\put(3,58){$C$}
\put(-269,117){$X$}
\put(-275,65){$Y$}
\put(-272,22){$Z$}
\put(-67,117){$X'$}
\put(-75,65){$Y'$}
\put(-69,22){$Z'$}
\caption{A different view of the delta move}
\label{fig:deltanu}
\end{figure}

We begin with the delta move as shown above and change it by an isotopy so that it now appears as pictured in Figure~\ref{fig:deltanu}.

\enlargethispage{\baselineskip}

%%% FIGURE 4: Delta move: nu %%%
%%%%%%%%%%%%%%%%%%%%

\noindent From this picture, we see that the delta move can be regarded as a pair of crossing changes of opposite sign which introduce a pair of small ``clasps'' $\put(0,-3){\includegraphics[height=10pt]{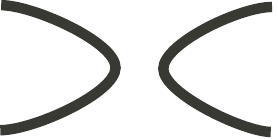}} \put(20,0){$\longrightarrow$} \put(40,-3){\includegraphics[height=10pt]{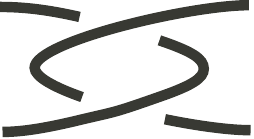}}$\hskip60pt.  

We may assume that the move takes place inside the large ball $B$ in $S^3$ of radius $\pi/2$ centered at $-1$.  If we think of $1$ and $-1$ as the north and south poles of $S^3$, then $B$ is just the southern hemisphere.  Figure~\ref{fig:deltanu} shows that inside $B$, the portions of $X$ and $X'$ lie along the great circle through $1$ and $i$ (so in fact $X$ and $X'$ coincide) while the portions of $Y$, $Z$, $Y'$ and $Z'$ lie close to the great circle $C$ through $j$ and $k$.

Outside $B$, the links $L$ and $L'$ coincide and, maintaining this coincidence, we move them into a more favorable position as follows.  First unknot $X=X'$ by a link homotopy, and move it to the rest of the great circle through $1$ and $i$.  Then by an isotopy move $Y$ and $Z$ into a small neighborhood of the great circle $C$, and position them so that their orthogonal projections to $C$ are Morse functions with just one critical point per critical value.  As intended, $Y'$ and $Z'$ move likewise outside $B$.  Note that each of the aforementioned clasps contributes two critical values to these projections, and we may assume that no other critical values fall between these two.

Now we use the asymmetric versions  $\g_L$  and $\g_{L'}$ of the characteristic maps $T^3 \to S^2$, identifying the target $S^2$ with the unit $2$-sphere in the purely imaginary $ijk$-subspace of the quaternions.  We interpret Figure~\ref{fig:deltanu} as showing the image of $S^3 - \{1\}$ under the stereographic projection $\pr_1$,  and view $i$ as the north pole of $S^2$.

\begin{figure}[h!]
	\vspace{-.1in}
	\includegraphics[height=130pt]{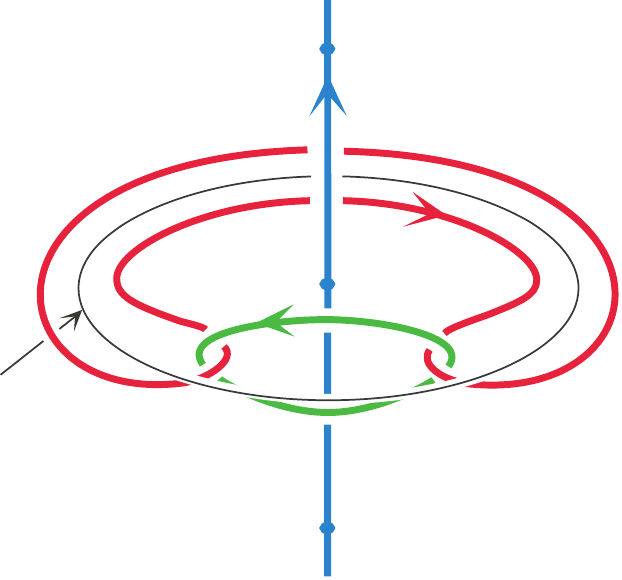}
	\put(-83,65){\small$-1$}
	\put(-75,118){\small$i$}
	\put(-82,11){\small$-i$}
	\put(-84,134){\small$\infty \leftrightarrow 1$}
	\put(-150,37){$C$}
	\put(-60,1){$X'$}
	\put(-88,27){$Z'$}
	\put(-12,36){$Y'$}
	\vspace{-.1in}
	\caption{Borromean Rings}
	\label{fig:borromean}
\end{figure}
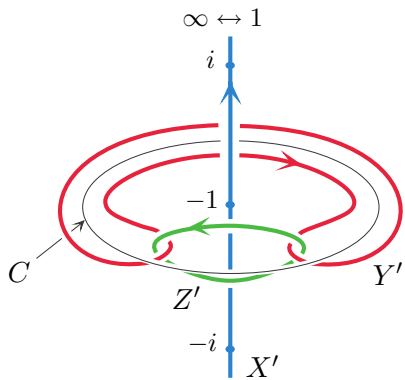

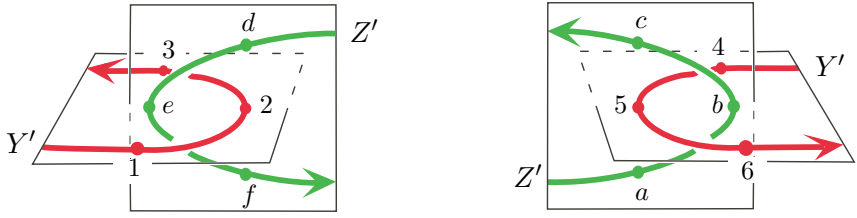
\begin{figure}[h!]
	\includegraphics[height=80pt]{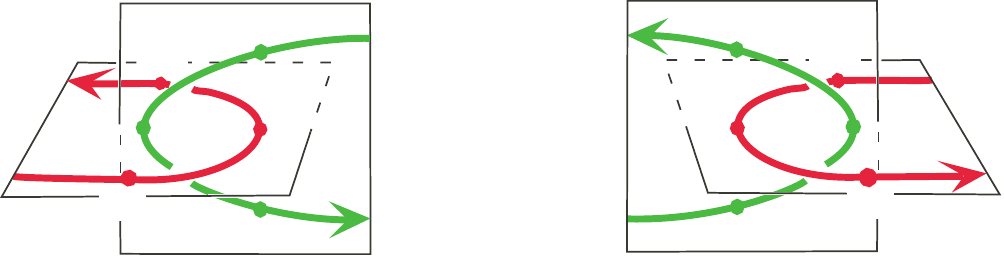}
	\put(-322,23){$Y'$}    
	\put(-16,52){$Y'$}    
	\put(-192,65){$Z'$}    
	\put(-130,10){$Z'$}    
	\put(-276,14){\small$1$}   
	\put(-226,38){\small$2$}    
	\put(-263,59){\small$3$}    
	\put(-55,60){\small$4$}    
	\put(-92,38){\small$5$}   
	\put(-44,13){\small$6$}    
	\put(-84,5){\small$a$}    
	\put(-55,38){\small$b$}    
	\put(-84,70){\small$c$}    
	\put(-233,69){\small$d$}    
	\put(-263,38){\small$e$}    
	\put(-233,4){\small$f$}    
	\caption{The two clasps}
	\label{fig:clasps}
\end{figure}

\begin{figure}[h!]
	\vspace{.15in}
\includegraphics[height=160pt]{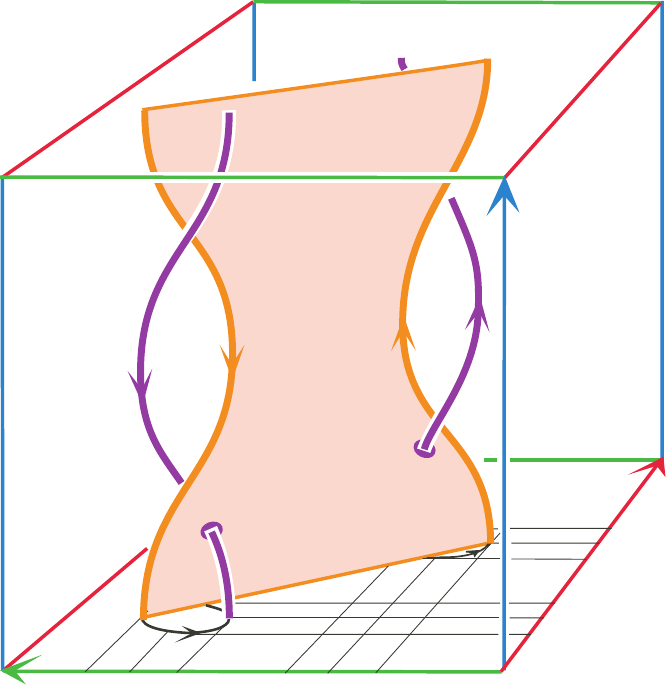}
\put(-175,90){$\tor$}
\put(-85,75){$\ba$}
\put(-142,-7){\small$f$}
\put(-131,-7){\small$e$}
\put(-120,-7){\small$d$}
\put(-94,-7){\small$c$}
\put(-83,-7){\small$b$}
\put(-72,-7){\small$a$}
\put(-8,34){\small$6$}
\put(-11,28){\small$5$}
\put(-15.5,24){\small$4$}
\put(-21,17){\small$3$}
\put(-26,11){\small$2$}
\put(-30,6){\small$1$}
\put(-40,127){$s$}
\put(2,57){$t$}
\put(-163,2){$u$}
\caption{Computing $e(\ba) = 2$}
\label{fig:annulus}
\end{figure}

It is straightforward to check that, under the genericity conditions imposed above, the point  $i \in S^2$  is a common regular value for $\g_L$  and $\g_{L'}$, and that the framed links $\bl$ and $\bl'$ in $T^3$ that are the inverse images of $i$ under these maps are for the most part the same.  In fact we will show that, up to isotopy
$
\bl' = \bl\cup\bl^*
$
where $\bl^*$ consists of a pair of oppositely oriented spiral perturbations of $S^1\times \text{pt} \times \text{pt}$, coming from the two clasps shown in Figure~\ref{fig:deltanu}, and that these two spirals bound an annulus $\ba$ in $T^3-\bl$ whose relative normal Euler class  $e(\ba)$ is $2$.  The argument will be given as we explain Figures \ref{fig:borromean}--\ref{fig:annulus}.

The discussion is independent of what the links $L$ and $L'$ look like outside the ball $B$ shown in Figure~\ref{fig:deltanu}.  The simplest case occurs when $L$ is the three-component unlink, and $L' = X' \cup Y' \cup Z'$  is the Borromean rings, whose stereographic image is shown in Figure~\ref{fig:borromean}.

%%% FIGURE 5: Borromean rings %%%
%%%%%%%%%%%%%%%%%%%%

\noindent We have, as before, that $X'$ is the great circle through $1$ and $i$, with image the blue vertical axis, while $Y'$ and $Z'$ lie in a small tubular neigborhood of the great circle $C$ through $j$ and $k$, with images shown in red and green.  In this circumstance, the general formula $\bl' = \bl \cup \bl^*$ has $\bl$ empty, and hence $\bl' = \bl^*$.

Figure~\ref{fig:clasps} shows enlargements of the two clasps between red $Y'$ and green $Z'$, with points on their segments labeled by numbers on $Y'$ and by letters on $Z'$, ordered consistently with their orientations.

%%% FIGURE 6: Clasps %%%
%%%%%%%%%%%%%%%%

Our job is to find the preimage $\bl'$ in $S^1 \times S^1 \times S^1$ of the regular value $i$ of the map $g_{L'}$,  which means we must find the points where the vector from the green $\pr_1(z\bx)$ to the red $\pr_1(y\bx)$ points straight up.

We pause to see the effect of right multiplication by $\bx$.  Let $x$ travel along $X = X'$ from $1$ to $-i$ to $-1$ to $i$ and back to $1$, which is the direction in which this component is oriented.  Then right multiplication by $\bx$ gradually rotates this component in the opposite direction.  In the image $3$-space, it looks like the vertical axis is moving downwards.  At the same time, the great circle $C$ through $j$ and $k$ is gradually rotated in the direction from $j$ towards $-k$.  A small tubular neighborhood of $C$ follows this rotation and twists as it goes, dragged by the downward motion of the vertical axis.

Now focusing on the left clasp, and starting with $x = 1$,  we see that the arrow $\vec{f2}$ from green $f$ to red $2$ points up.  As $x$ moves around $X'$ from $1$ towards $-i$ (up on the blue vertical axis) a loop of upward pointing vectors is traced out, passing successively through $\vec{e1}$,  $\vec{d2}$, $\vec{e3}$, and finally back to $\vec{f2}$.  In Figure~\ref{fig:annulus}, the $3$-torus is depicted as a cube, in which this loop is shown near the front left corner of the bottom red-green square, traced in a counterclockwise direction.  When the progression of $x(s)$ values is taken into account, we get the orange spiral curve shown above this loop.  This is one component of $\bl'$, and is oriented according to the convention for framed links.

%%% FIGURE 7: Annulus %%%
%%%%%%%%%%%%%%%%%

Focusing on the right clasp and repeating the above procedure, we get the orange spiral curve shown at the right rear, lying over the loop $\vec{a5}\,\vec{b6}\,\vec{c5}\,\vec{b4}$.  This is the other component of $\bl'$, and is also oriented according to convention.\foot{In fact, one need not go through the careful analysis to determine the orientations of the spirals.  All that is important is that they are oppositely oriented, which follows from the fact that the pairwise linking numbers for the Borromean rings are zero.}  

Together, the components of $\bl'$ bound the orange annulus $\ba$ shown in the cube.
This annulus is constructed as follows.  Each point on the left spiral loop is joined  to the point on the right spiral loop at the same height $s$ by a straight line segment in the 2-torus $s \times S^1 \times S^1$.  The $t$-coordinate moves steadily so that $y(t)$ travels along the arc of $Y'$ which lies in the ball $B$ (joining $3$  to $4$) and the $u$-coordinate moves steadily so that $z(u)$ travels along the arc of $Z'$ which lies in $B$ (joining $d$ to $c$).  It is easy to see that this annulus $\ba$ is embedded in the $3$-torus and that, even in the general case where $\bl$ is not empty, it would still be disjoint from $\bl$.

Now since $\ba$ lies in $\tor$, its relative normal Euler class $e(\ba)$ (when viewed as a surface in $\tor\times[0,1]$) can be computed as the intersection number of $\ba$ with the inverse image of any other regular value of $\g_{L'}$.  In particular, the point $-i$ is also a regular value, and its inverse image $\widehat\bl'$ is calculated from an analysis of the clasps, just as we did for the inverse image of $i$.  It consists of two spirals, which can be obtained from the spirals in $\bl'$ by moving them half-way along in the vertical (blue) direction.  We show $\widehat\bl'$ in Figure~\ref{fig:annulus} as a pair of purple spirals, which are oriented the same way as the orange spirals. 

\enlargethispage{\baselineskip}

It is seen in this figure that the purple spirals pierce the orange annulus twice in the positive direction, confirming that 
\[
e(\ba)  \ =  \ \ba \,\dt\, \widehat\bl'  \ =  \ 2,
\]
and completing the sketch of the topological proof of \thmref{thm:A}.

\clearpage

%%%%%%%%%%%%%%%%%%%%%%%%%
%% SECTION 5: Sketch of the algebraic proof of Theorem A
%%%%%%%%%%%%%%%%%%%%%%%%%

%%%%%%%%%%%%%%%%%%%%%%%%%%%%%%
\section{Sketch of the algebraic proof of \thmref{thm:A}}\label{sec:alg_proof}
%%%%%%%%%%%%%%%%%%%%%%%%%%%%%%

This proof is organized around the following key diagram:
\begin{equation}\label{eqn:keydiagram}
	\begin{diagram}[size=2em]
		[[S^1 \cup S^1 \cup S^1, S^3]]  & \lTo & \calh(3)  \\
		\dTo^g & & \dTo_{\calg} \\
		[S^1 \times S^1 \times S^1, S^2] & \lTo & \pi_1\M \\
	\end{diagram}
	\tag{$*$}
\end{equation}
The left half represents the geometric-topological problem we are trying to solve, and is devoid of algebraic structure.  The right half represents the algebraic structures that we impose on the left half via the two horizontal maps in order to solve the problem.

In the upper left corner of \eqref{eqn:keydiagram}  we have the set of link-homotopy classes of three-component links in the 3-sphere $S^3$, and in the lower left corner the set of homotopy classes of maps of the 3-torus to the 2-sphere.  The vertical map $g$ between them assigns to the link-homotopy class of $L$ the homotopy class of its characteristic map $g_L$.  \thmref{thm:A} describes $g$ and asserts that it is one-to-one.

In the upper right corner of \eqref{eqn:keydiagram} 
we have the group $\calh(3)$ of link-homotopy classes of three-component string links.  A $k$-component {\itb string link} consists of $k$ oriented intervals embedded disjointly in a cube, with their tails on the bottom face, their tips on the top face directly above their tails, and their interiors in the interior of the cube. The terminology was coined by Habegger and Lin~\cite{HabeggerLin}.  The product of two $k$-component string links with endpoints in a common position is given by stacking the second one on top of the first.  When a string link moves by a link homotopy, each strand is allowed to cross itself, while different strands must remain disjoint, just as for links.  Then the above product induces a group structure on the set $\calh(k)$ of link homotopy classes of $k$-component string links.

%%% FIGURE ?: H(3) Generators %%%
\begin{figure}[h!]
\includegraphics[height=150pt]{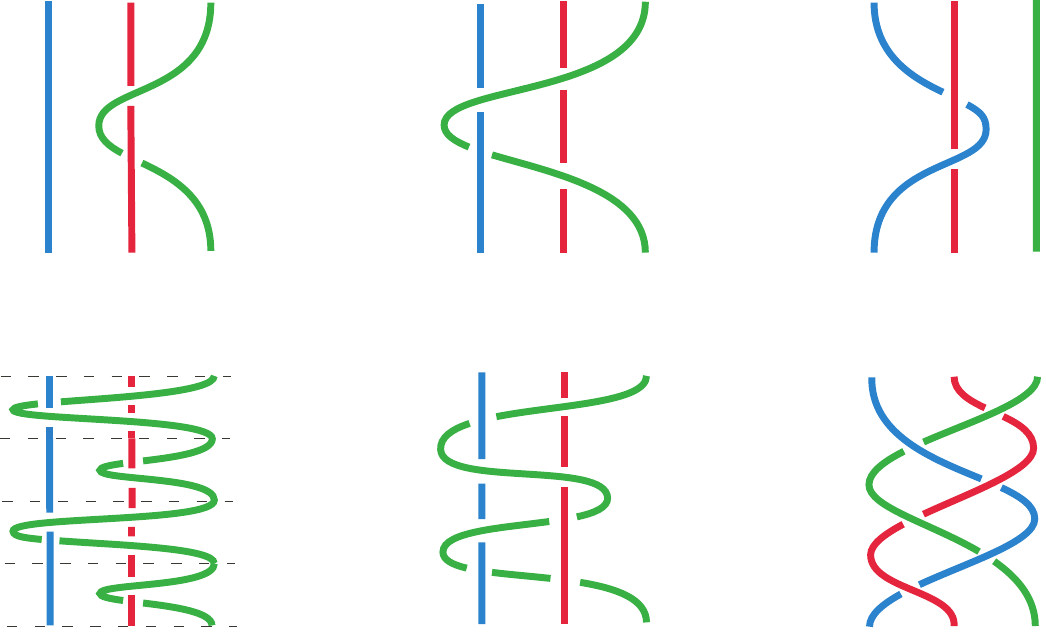} % replaces old FigGenerators.pdf
\put(-225,77){\small$P$}
\put(-120,77){\small$Q$}
\put(-25,77){\small$R$}
\put(-80,29){$=$}
\put(-175,29){$=$}
\put(-165,-15){$PQP^{-1}Q^{-1} \ = \ B$}
\put(-263,7){\footnotesize $P$}
\put(-263,22){\footnotesize $Q$}
\put(-265,37){\footnotesize $P^{-1}$}
\put(-265,52){\footnotesize $Q^{-1}$}
\caption{Generators for $\calh(3)$}
\label{fig:h3_generators}
\end{figure}
%%%%%%%%%%%%%%%%%

Following Habegger and Lin, we derive the following presentation for $\calh(3)$:
\[
	\begin{aligned}
	\calh(3) \ = \ \langle P, Q, R, B \ | \ &[P,Q] = [Q,R] = [R,P] = B, \\
	&[P,B] = [Q,B] = [R,B] = 1 \rangle
	\end{aligned}
	%\calh(3) \ = \ \langle P, Q, R, B \ | \ [P,Q] = [Q,R] = [R,P] = B, \,[P,B] = [Q,B] = [R,B] = 1 \rangle
\]
The string links $P$, $Q$, $R$ and $B$ are those shown in Figure~\ref{fig:h3_generators}.

Using this presentation, elements of $\calh(3)$ can be written uniquely in the form
$$
P^pQ^qR^rB^\mu  \ , \quad \text{for} \ p,q,r,\mu\in\bz.
$$
Two elements $P^pQ^qR^rB^\mu$ and $P^{p^*}Q^{q^*}R^{r^*}B^{\mu^*}$ are conjugate if and only if $p = p^*$, $q = q^*$, $r = r^*$ and $\mu \equiv \mu^*$ mod $\gcd(p,q,r)$.   

\begin{figure}[h!]
\includegraphics[height=90pt]{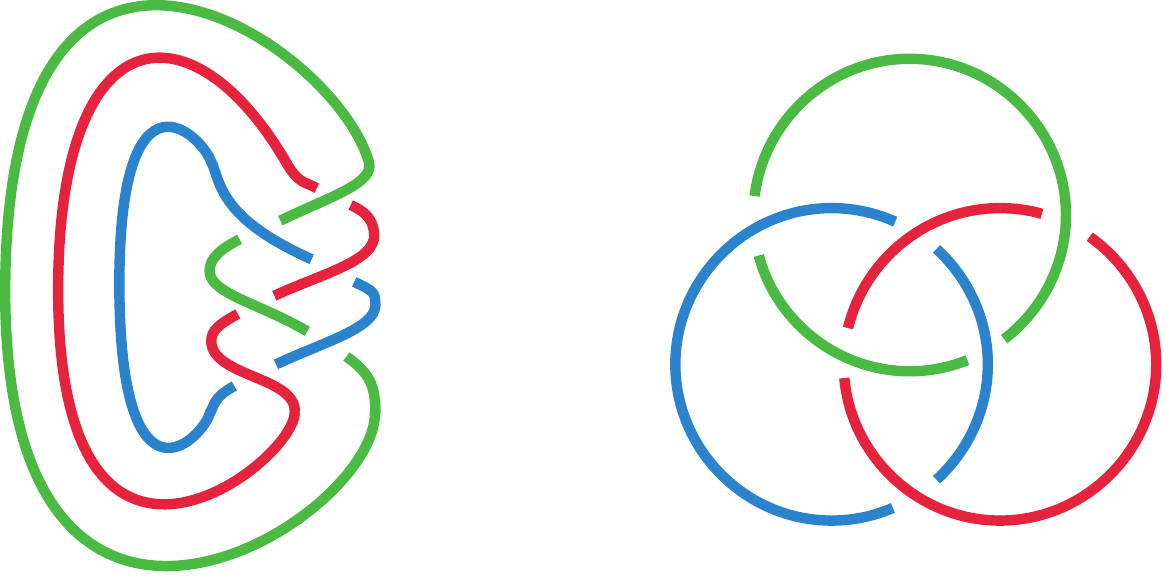} 
\put(-103,40){$=$}\put(-125,7){\small$\widehat B$}
\caption{The string link $B$ closes up to the Borromean rings $\widehat B$}
\label{fig:B_closure}
\end{figure}

A string link $S$ can be closed up to a link $\widehat{S}$  by joining the tops of the strands to their bottoms outside the cube, without introducing any more crossings.  For example, the closure of the three-component string link $B$ is the Borromean rings, as shown in Figure~\ref{fig:B_closure}.

%%% FIGURE ?: B closure %%%
%%%%%%%%%%%%%%%%%

Thus the Borromean rings, a ``primitive example'' in the world of links, is the closure of a string link which is itself a commutator of simpler string links.

The closing-up operation descends to link homotopy classes and provides the upper horizontal map in \eqref{eqn:keydiagram}.    It is easily seen that the closure of the string link $P^pQ^qR^rB^\mu$ has pairwise linking numbers $p$, $q$ and $r$ and Milnor invariant $\mu\mod{\gcd(p,q,r)}$.  It follows that this closing-up map is onto.  Furthermore, its point inverse images are the conjugacy classes in $\calh(3)$, a special circumstance for links with three components which fails for four or more components.

In the lower right corner of \eqref{eqn:keydiagram} we have the union of the fundamental groups of the components of the space of continuous maps of the $2$-torus to the $2$-sphere, with one group $\pi_1 \text{Maps}_p(S^1 \times S^1, S^2)$ for each choice of degree $p$ of these maps.  The work of Fox~\cite{Fox} on torus homotopy groups (the case $p = 0$), and its generalization by Larmore and Thomas~\cite{LarmoreThomas}, provides explicit presentations for these groups:
\[
	\begin{aligned}
	\pi_1\Mp \ = \ \langle U_p, V_p, W_p \st &[U_p,V_p] = W_p^2\,,\,W_p^{2p} = 1, \\
	&[U_p,W_p] = [V_p,W_p] = 1 \rangle
	\end{aligned}
\]

Using this presentation, elements of $\pi_1\Mp$ can be written uniquely in the form
\[
U_p^q\, V_p^r\, W_p^\nu  \ , \quad \text{for} \ q,r,\nu\in\bz \text{ with } 0 \leq \nu < 2|p|.
\]
Two elements $U_p^q \, V_p^r \, W_p^\nu$ and $U_p^{q^*} V_p^{r^*} W_p^{\nu^*}$ are conjugate if and only if $q = q^*$, $r = r^*$ and $\nu \equiv \nu^*$ mod $2\gcd(p,q,r)$.

A direct argument using framed links shows how the generators $U_p$, $V_p$ and $W_p$ of $\pi_1\Mp$ can be represented by specific maps $S^1 \times S^1 \times S^1 \to S^2$,  all agreeing with some fixed map of degree $p$ on $* \times S^1 \times S^1$.  In addition to this common feature, the representatives for $U_p$, $V_p$, $W_p$ have degrees $q = 1, 0, 0$ on $S^1 \times * \times S^1$, degrees $r = 0, 1, 0$ on $S^1 \times S^1 \times *$,  and Pontryagin invariants $\nu = 0, 0, 1$ relative to the chosen basepoint for $\pi_1\Mp$.

The lower horizontal map in the key diagram takes a homotopy class of based loops in the space $\M$,  ignores basepoints, identifies  $S^1 \times S^1$ with $* \times S^1 \times S^1$,  and then regards this class as a homotopy class of maps of $S^1 \times S^1 \times S^1 \to S^2$ in the usual way.  In particular, the element $U_p^q \, V_p^r \, W_p^\nu$ of  $\pi_1\Mp$ is taken to a map with degrees $p$, $q$ and $r$ on the $2$-dimensional coordinate subtori, and -- {\sl this is the key observation} -- with Pontryagin invariant  $\nu$ mod $2\gcd(p,q,r)$ relative to the image of the basepoint map.  This lower horizontal map is onto, and point inverse images are conjugacy classes in $\pi_1 \Mp$.

The final step will be to define the vertical map $\calg$ on the right side of the key diagram to make the whole diagram commutative, and to be a group homomorphism, insofar as possible.  The hedge ``insofar as possible'' refers to the fact that we have a group $\calh(3)$ in the upper right corner of the diagram, but only a union of groups $\pi_1 \text{Maps}_p(S^1 \times S^1, S^2)$ in the lower right corner.  We deal with this disparity by demoting $\calh(3)$ to a union of groups as follows.

Let $\calh_0(3)$ denote the subgroup of $\calh(3)$ consisting of string links with the second and third strands unlinked, and with presentation:
\[
\calh_0(3) = \langle Q,R,B \ | \ [Q,R] = B\,,\,[Q,B] = [R,B] = 1 \rangle.
\]
Consider the left cosets $\calh_p(3) = P^p \calh_0(3)$  of $\calh_0(3)$,  and convert each of them into a subgroup isomorphic to $\calh_0(3)$ by using left translation to transfer the group structure from subgroup to coset.  Adopting the notations $Q_p = P^p Q$, $R_p = P^p R$ and $B_p = P^p B$  for the generators of $\calh_p(3)$ in this borrowed group structure, we get the presentation
\[
\calh_p(3) = \langle Q_p, R_p, B_p \ | \ [Q_p,R_p] = B_p\,,\,[Q_p,B_p] = [R_p,B_p] = 1 \rangle.
\]

We are now ready to define the vertical map $\calg$ on the right side of \eqref{eqn:keydiagram} so as to make the whole diagram commutative, and at the same time to be a union of homomorphisms from the groups $\calh_p(3)$ to the groups $\pi_1 \text{Maps}_p(S^1 \times S^1, S^2)$.

\pagebreak

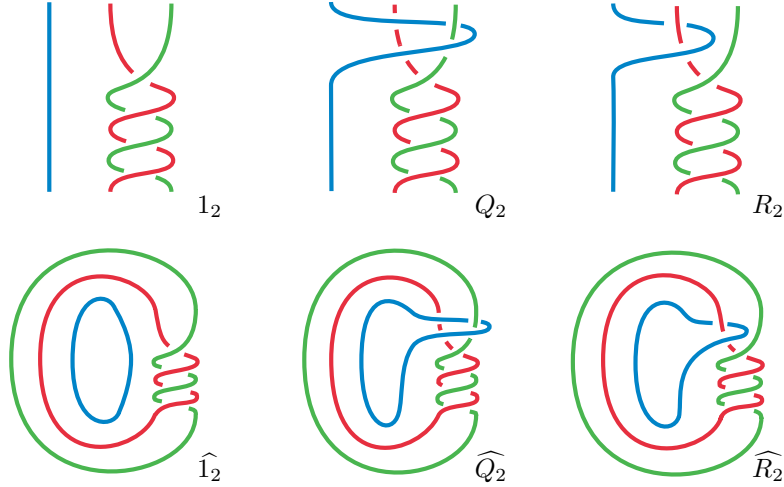
\begin{figure}[h!]
\includegraphics[height=180pt]{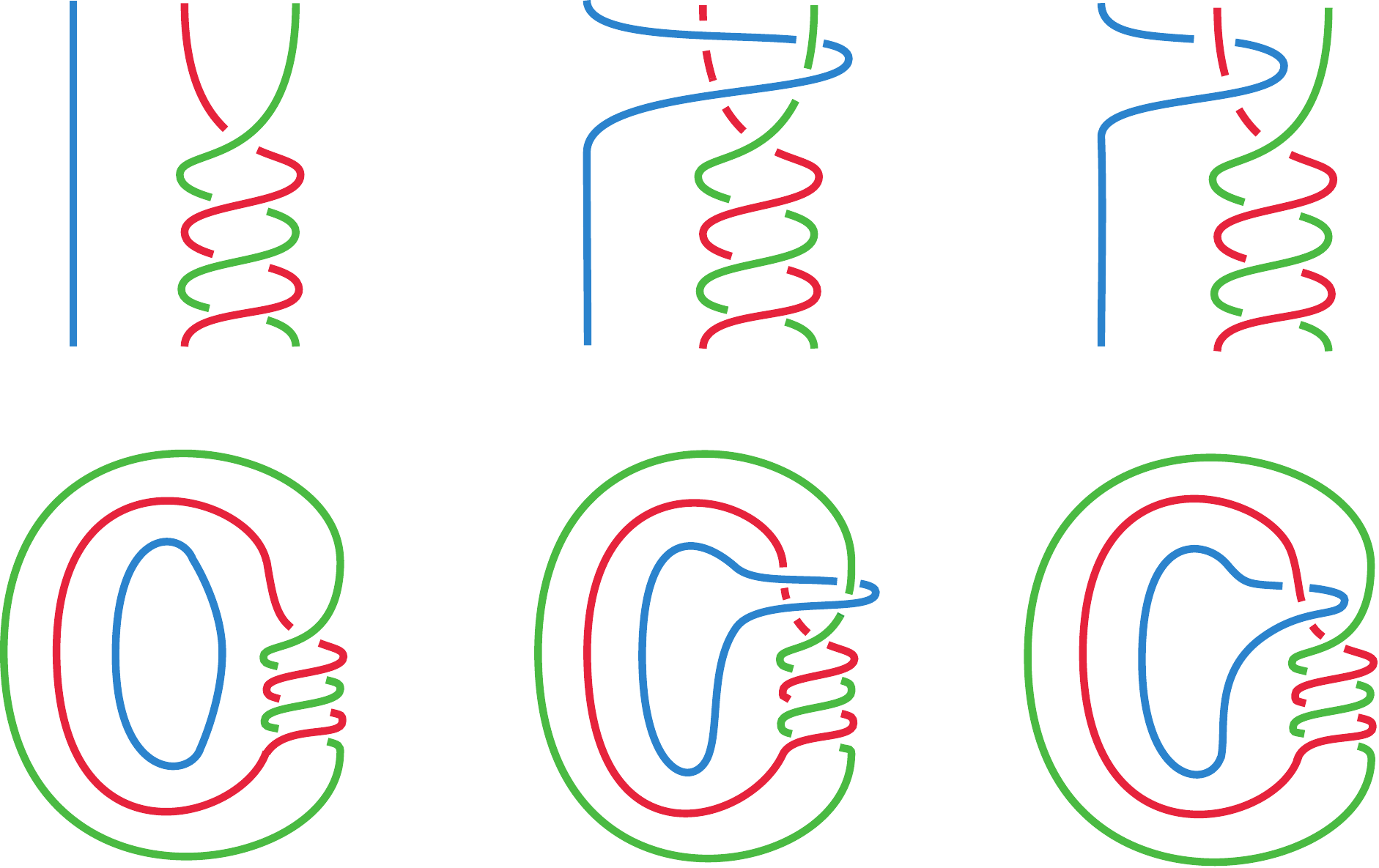} 
\put(-215,100){\small$1_2$}
\put(-110,100){\small$Q_2$}
\put(-5,100){\small$R_2$}
\put(-215,-0){\small$\widehat{1_2}$}
\put(-110,0){\small$\widehat{Q_2}$}
\put(-5,0){\small$\widehat{R_2}$}
\caption{Closing up the generators of $\calh_2(3)$}
\label{fig:string_closed}
\end{figure}

To do this, we start with specific string links to represent the elements of $\calh_p(3)$.  For purposes of illustration, we choose $p = 2$,  and show in Figure~\ref{fig:string_closed} the string  links $1_2 = P^2$, $Q_2 = P^2 Q$ and $R_2 = P^2 R$,  and under them the three-component links we get by closing them up.

%%% FIGURE ?: Generators for H_p(3) %%%
%%%%%%%%%%%%%%%%%

The corresponding characteristic maps from $S^1 \times S^1 \times S^1 \to S^2$ all restrict to the same map of degree $2$ on $* \times S^1 \times S^1$, and therefore all three represent elements of the fundamental group $\pi_1\text{Maps}_2(S^1 \times S^1, S^2)$ based at this map.  We denote these three images by  $\calg_2(P^2)$, $\calg_2(P^2 Q)$ and $\calg_2(P^2 R)$,  with the intent of forcing commutativity in the key diagram.  In fact, we can do this for all the string links $P^2 Q^q R^r B^\mu$, and a simple geometric argument shows that composition in the group $\calh_2(3)$ carries over in this way to multiplication in the group $\pi_1 \text{Maps}_2(S^1 \times S^1, S^2)$.  

Furthermore, a direct argument using framed links shows that the elements $\calg_2(P^2 Q)$ and $\calg_2(P^2 R)$ may serve as the elements $U_2$ and $V_2$ in the above presentation for $\pi_1\text{Maps}_2(S^1 \times S^1, S^2)$,  so that $\calg_2$ takes $P^2 Q$ to $U_2$ and $P^2 R$ to $V_2$.  It then follows that $\calg_2$ maps $P^2 B = [P^2Q, P^2R]$ to $[U_2,V_2] = W_2^2$.  

\enlargethispage{\baselineskip}

The value $p = 2$ used above was just for purposes of illustration, and the corresponding results are true for all values of $p$.  Thus we have defined the vertical map $\calg$ on the right side of our key diagram to be a union of homomorphisms $\calg_p :  \calh_p(3) \to \pi_1\Mp$ making the whole diagram commutative.  

Now let $L$ be any three-component link in $S^3$ with pairwise linking numbers $p$, $q$ and $r$ and Milnor invariant $\mu$.  Then $L$ is link homotopic to the closure of  $P^p Q^q R^r B^\mu$.  By commutativity of the key diagram, the homotopy class of the characteristic map $g_L$ is the image under the lower horizontal map of the element $\calg_p(P^p Q^q R^r B^\mu) = U_p^q \, V_p^r \, W_p^{2\mu}$  of $\pi_1\Mp$,  and therefore has 
Pontryagin invariant $2 \mu$, as desired.

This completes our sketch of the algebraic proof of \thmref{thm:A}.
\pagebreak

%%%%%%%%%%%%%%%%%%%%%%%%%
%% SECTION 6: Sketch of the proof of Theorem B
%%%%%%%%%%%%%%%%%%%%%%%%%

\section{Sketch of the proof of Theorem~\ref{thm:B}}\label{sec:thmB}
Let $L$ be a 3-component link in $S^3$ with pairwise linking numbers $p$, $q$ and $r$ all zero.  We saw in Theorem~\ref{thm:A} that these numbers are the degrees of the characteristic map $g_L: T^3 \to S^2$ on the 2-dimensional coordinate subtori.  Thus $g_L$ is homotopic to a map which collapses the 2-skeleton of $T^3$, and so is in effect a map of $S^3 \to S^2$.  The Hopf invariant of this map, which we will regard as the Hopf invariant of $g_L$, is equal to Pontryagin's $\nu$-invariant comparing $g_L$ to the constant map, and we will denote this by $\nu(g_L)$.

To calculate this Hopf invariant, we adapt J.\,H.\,C.~Whitehead's integral formula for the Hopf invariant of a map from $S^3 \to S^2$ to the case of a map from $T^3 \to S^2$, and show how to make the calculation explicit.

Whitehead~\cite{Whitehead} expressed the Hopf invariant of a map $f: S^3 \to S^2$ as follows.  Let $\omega$ be the area 2-form on $S^2$, normalized so that $\int_{S^2} \omega = 1$.  Then its pullback $f^* \omega$ is a closed 2-form on $S^3$ which is exact because $H^2(S^3; \br) = 0$.  Hence $f^*\omega = d\alpha$ for some 1-form $\alpha$ on $S^3$, and Whitehead showed that
\[
	\text{Hopf}(f) \ = \ \int_{S^3} \alpha \wedge f^*\omega,
\]
the integral being independent of the choice of $\alpha$.

We recast Whitehead's formula in terms of vector fields by letting $V_f$ be the vector field on $S^3$ corresponding in the usual way to the 2-form $f^* \omega$.  Then $V_f$ is divergence-free, since $f^* \omega$ is closed, and is in fact in the image of curl since $f^* \omega$ is exact.  Thus $V_f = \nabla \times W$ for some vector field $W$ on $S^3$, and the integral formula for the Hopf invariant becomes
$$
\text{Hopf}(f) \ = \ \int_{S^3} W \dtw V_f \ d(\text{vol}),
$$
independent of the choice of $W$.  

To make Whitehead's formula more explicit, one needs a way to produce a vector field $W$ whose curl is $V_f$.  On $\br^3$ this can be done by viewing $V_f$ as a flow of electric current and then calculating the corresponding magnetic field $\text{BS}(V_f)$ using the classical formula of Biot and Savart~\cite{BiotSavart}:
\[
	\text{BS}(V_f) \ = \ -\nabla \times \text{Gr}(V_f),
\] 
where $\text{Gr}$ is the Green's operator that inverts the vector Laplacian.  Then by Amp\`ere's Law, \mbox{$\nabla \times \text{BS}(V_f) = V_f$}.

The explicit formula for the Green's operator is well known on $\br^3$, namely convolution with the fundamental solution $\varphi(r) = -1/(4\pi r)$ of the scalar Laplacian, and hence
\[
	\text{BS}(V_f)(y) \ = \ \int_{\br^3} V_f(x) \times \nabla_y\varphi\left(\|y-x\|\right)\ dx,
\]
assuming that $V$ is compactly supported, to guarantee that the integral converges.   The corresponding formula on $S^3$ was given by DeTurck and Gluck~\cite{DeTurckGluck} and by Kuperberg~\cite{Kuperberg}.

We can summarize the calculation of the Hopf invariant in the single formula
\[
	\text{Hopf}(f) \ = \ \int_{S^3} \text{BS}(V_f) \dtw V_f \ d(\text{vol})
\]
which was Woltjer's original expression for the helicity of the vector field $V_f$.

A routine check shows that the above formula, with the integration over $T^3$ instead of $S^3$, yields the value of the Hopf invariant of the characteristic map $g_L: \tor \to S^2$.  This provides a formula for Pontryagin's $\nu$-invariant, and portrays it as the helicity of the associated vector field $V_{g_L}$, which for simplicity we denote by $V_L$:
\begin{equation}\label{eqn:PontryaginHelicity}
	\nu(g_L) = \text{Hopf}(g_L) = \text{Hel}(V_L) = \int_{T^3} \text{BS}(V_L) \dtw V_L \ d(\text{vol}).
	\tag{$\ddag$}
\end{equation}
A straightforward calculation shows that
\[
	V_L \ = \ \frac{F_t \times F_u \dtw F}{4\pi \|F\|^3} \, {\partial}/{\partial s} + \frac{F_u \times F_s \dtw F}{4\pi \|F\|^3} \, {\partial}/{\partial t} + \frac{F_s \times F_t \dtw F}{4\pi \|F\|^3} \, {\partial}/{\partial u},
\]
where $F: T^3 \to \br^3 - \{0\}$ is the map defined in Section~\ref{sec:characteristicmap}, and where subscripts denote partial derivatives.

Therefore, to make the integral formula (\ref{eqn:PontryaginHelicity}) explicit, it remains to obtain an explicit formula for the Biot-Savart operator on the 3-torus.  As in the case of $\br^3$, this depends on having an explicit formula for the fundamental solution of the scalar Laplacian.

%%%%%%%%%%%%%%%%%%%%%%
\part{The fundamental solution of the Laplacian} 
\label{sub:fund_soln}
%%%%%%%%%%%%%%%%%%%%%%

\begin{proposition}\label{prop:fund_soln}
The fundamental solution of the scalar Laplacian on the 3-torus $\tor$ is given by the formula
	\[
		\varphi(x,y,z) \ = \ -\,\frac{1}{8\pi^3}\,\sum \, \frac{e^{i(\ell x+my+nz)}}{\ell^2+m^2+n^2}
	\]
where the sum is over all integers $\ell$, $m$ and $n$ $($positive and negative$)$ with the exception of $\ell=m=n=0$.
\end{proposition}

%-\frac{1}{8\pi^3}

Even though we have expressed $\varphi$ in terms of complex exponentials, the value of $\varphi$ is real for real values of $x$, $y$ and $z$ because of the symmetry of the coefficients, and can therefore also be expressed as a Fourier cosine series.

% \vskip -.2in

%%% FIGURE 6: Fundamental Solution of Laplacian %%%
\begin{figure}[h!]
\includegraphics[height=200pt]{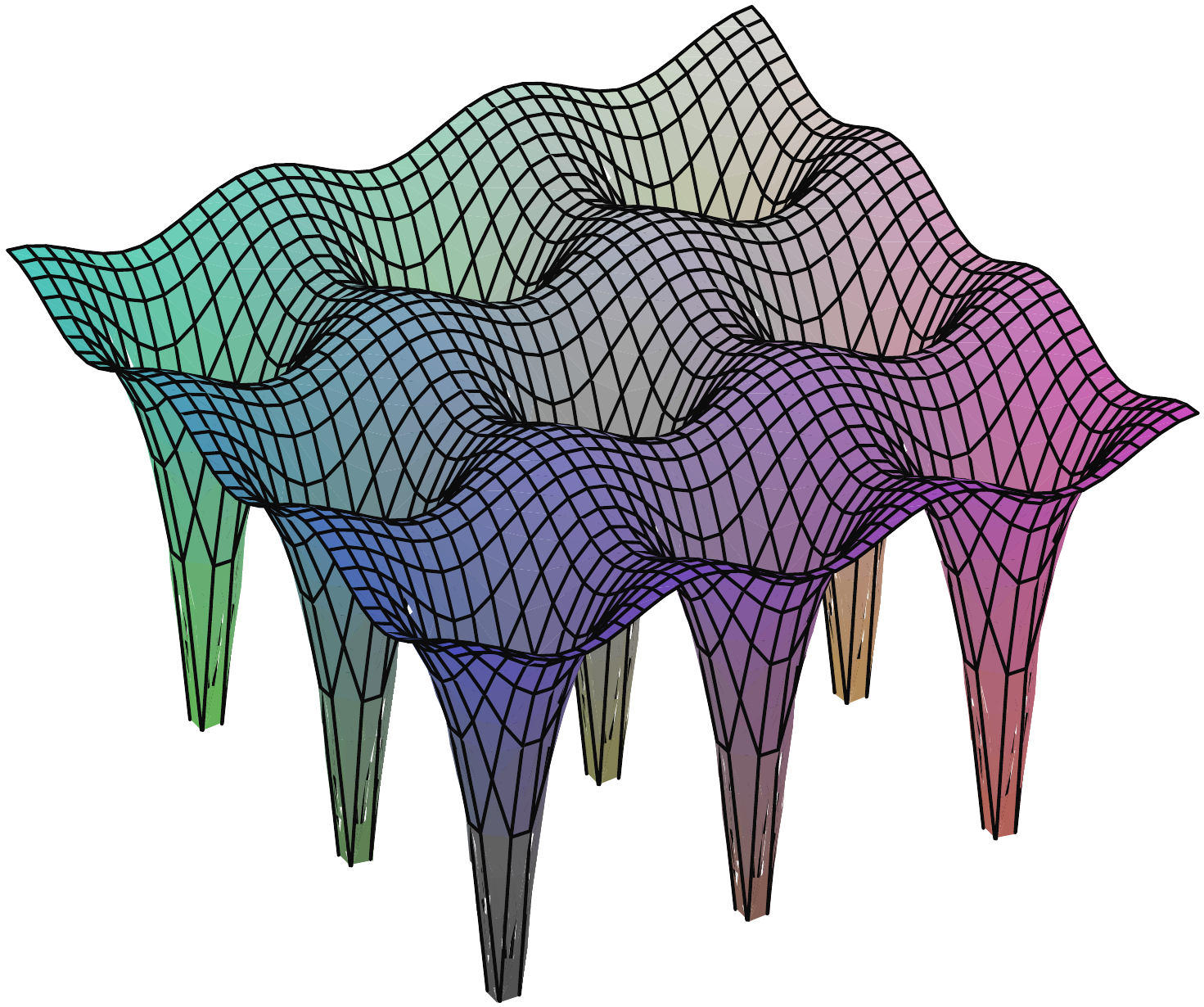}
% \vskip -.2in
\caption{Fundamental solution of the scalar Laplacian on the $2$-torus}
\label{fig:fund_sol}
\end{figure}
%%%%%%%%%%%%%%%%%

Figure~\ref{fig:fund_sol} shows the graph of the corresponding fundamental solution of the scalar Laplacian on the 2-torus $T^2$, displayed over the range $-3\pi \leq x, y \leq 3\pi$.  If we think of the 2-torus as obtained from a square by identifying opposite sides, then this shows the function $\varphi$ to have a negative infinite minimum at the single vertex, two saddle points in the middle of the two edges, and a maximum in the middle of the square.  Presumably the fundamental solution $\varphi$ on the 3-torus displays a corresponding Morse-like behavior.

%%%%%%%%%%%%%%%%%%%%%%
\part{Completing the proof of Theorem~\ref{thm:B}} 
\label{sub:ThmBCompleted}
%%%%%%%%%%%%%%%%%%%%%%

Now we have an explicit formula for the fundamental solution $\varphi$ of the scalar Laplacian on the 3-torus $T^3$, and we know that both the scalar and vector Green's operators act by convolution with $\varphi$.  In particular, if $V$ is a smooth vector field on $T^3$, then $\text{Gr}(V) = V * \varphi$, that is,
\[
	\text{Gr}(V)(\tau) \ = \ \int_{T^3} V(\sigma) \varphi(\tau - \sigma) \ d\sigma.
\]

To obtain the formula for the magnetic field $\text{BS}(V)$, we take the negative curl of the above formula and get
\begin{align*}
	\text{BS}(V)(\tau) \ = \ -\nabla_\tau \times \text{Gr}(V)(\tau) \ & = \ -\int_{T^3} \nabla_\tau \times \left(V(\sigma) \varphi(\tau - \sigma)\right)d\sigma \\
	& = \ \int_{T^3} V(\sigma) \times \nabla_\tau \varphi(\tau - \sigma) \ d\sigma.
\end{align*}
Then the helicity of $V$ is given by
\begin{align*}
	\text{Hel}& (V) \ = \ \int_{T^3} V(\tau) \dtw \text{BS}(V)(\tau) \ d\tau \\
	& = \ \int_{T^3 \times T^3} V(\sigma) \times V(\tau) \dtw \nabla_\sigma \varphi(\sigma - \tau) \ d\sigma \, d\tau.
\end{align*}

Applying this to the vector field $V_L$ associated with our 3-component link $L$, we get the desired formula for the Pontryagin invariant $\nu$ of $g_L$:
\begin{align*}
	\nu(g_L) \ &= \ \text{Hopf}(g_L) = \text{Hel}(V_L) = \int_{T^3} \text{BS}(V_L) \dtw V_L \ d(\text{vol}) \\
	&= \int_{T^3 \times T^3} V_L(\sigma) \times V_L(\tau) \dtw \nabla_\sigma \varphi(\sigma-\tau) \ d \sigma \, d \tau.
\end{align*}

Hence by Theorem~\ref{thm:A}, Milnor's $\mu$-invariant of the 3-component link $L$ is given by 
\[
	\mu(L) \ = \ \frac12\,\nu(g_L) \ = \ \frac12\, \int_{T^3 \times T^3} V_L(\sigma) \times V_L(\tau) \dtw \nabla_\sigma \varphi(\sigma-\tau) \ d\sigma \, d\tau,
\]
completing the proof of \thmref{thm:B}.
% subsection ThmBCompleted (end)

%%%%%%%%%%%%%%%%%
%% REFERENCES
%%%%%%%%%%%%%%%%%

\providecommand{\bysame}{\leavevmode\hbox to3em{\hrulefill}\thinspace}
\providecommand{\MR}{\relax\ifhmode\unskip\space\fi MR }
% \MRhref is called by the amsart/book/proc definition of \MR.
\providecommand{\MRhref}[2]{%
  \href{http://www.ams.org/mathscinet-getitem?mr=#1}{#2}
}
\providecommand{\href}[2]{#2}

\bigskip

\noindent\bf\tt
deturck@math.upenn.edu \\
gluck@math.upenn.edu \\
rako@math.upenn.edu \\
pmelvin@brynmawr.edu \\
shonkwil@math.upenn.edu \\
dvick@math.upenn.edu

\end{document}